\documentclass{amsart}
\usepackage[ansinew]{inputenc}
\usepackage[all]{xy}
\usepackage{amsmath,latexsym,amssymb,verbatim,amsthm}

\input arrow.tex

\newcommand{\N}{\mathbb{N}}

\newtheorem{lema}{Lemma}[section]
\newtheorem{teo}[lema]{Theorem}
\newtheorem{pro}[lema]{Proposition}
\newtheorem{defi}[lema]{Definition}

\newtheorem{com}[lema]{Remark}

\newtheorem{den}[lema]{Notation}
\newtheorem{teo*}{Theorem}

\begin{document}

\title[Invariant  Integral and  Algebraic Harmonic Analysis]{Invariant  Integral on Classical Groups and  Algebraic Harmonic Analysis}

\author{Amelia \'{A}lvarez}
\address{Departamento de Matem\'{a}ticas, Universidad de Extremadura,
Avenida de Elvas s/n, 06071 Badajoz, Spain}
\email{aalarma@unex.es}

\author{Carlos Sancho}
\address{Departamento de Matem\'{a}ticas, Universidad de Salamanca,
Plaza de la Merced 1-4, 37008 Salamanca, Spain}
\email{mplu@usal.es}

\author{Pedro Sancho}
\address{Departamento de Matem\'{a}ticas, Universidad de Extremadura,
Avenida de Elvas s/n, 06071 Badajoz, Spain} \email{sancho@unex.es}


\begin{abstract}
Let $G= {\rm Spec}\, A$ be a linearly reductive group and let $w_G \in A^*$ be the invariant integral on $G$. We establish the algebraic harmonic analysis on $G$ and we compute $w_G$ when $G= Sl_n, Gl_n, O_n, Sp_{2n}$ by geometric arguments and by means of the Fourier transform.
\end{abstract}

\maketitle

\section*{Introduction}

An affine $k$-group $G = {\rm Spec}\, A$ is linearly semisimple (that is, linearly reductive) if and only if ${A}^*$ splits into the form ${A}^* = {k} \times {B}^*$ as $k$-algebras, where the first projection $\pi_1\colon {A}^*\to {k}$ is the morphism $\pi_1(w):=w(1)$ (\cite[2.6]{Reynolds}). The linear form $w_G:=(1,0) \in {k} \times {B}^*=A^*$ will be referred
to as the {\it invariant integral} on $G$.

The aim of this article is to establish the algebraic harmonic analysis on $G$ and the explicit calculation of $w_G$ when $G=Sl_n, Gl_n$, $O_n$, $Sp_{2n} $  (${\rm char}\, k=0$) by geometric arguments and by means of the Fourier transform, which is defined below.

An affine $k$-group $G = {\rm Spec}\, A$ is linearly semisimple if and only if ${A}^* = \prod_i {A}^*_i$, where $A_i^*$ are finite simple $k$-algebras (\cite[6.8]{Amel}). On $\tilde A:=\oplus_iA_i^*\subseteq A^*$ one has the non-singular trace metric and its associated polarity $\phi\colon \tilde A \simeq A$. This morphism is essentially characterized by being a morphism of left and right $A^*$-modules or equivalently of left and right $G$-modules. Let  $* : A \to A$, $a \mapsto a^*$ be the morphism induced by the morphism $G \to G$, $g \mapsto g^{-1}$. We prove that the Fourier Transform $$F\colon {A} \to {A}^*,\, F(a)= w_G(a^*\cdot -),$$ where $w_G(a^* \cdot -)(b):=w_G(a^* \cdot b)$, is the inverse morphism of $\phi$, because is a morphism of left and right $G$-modules.  The product operation on ${A}^*$ defines, via the Fourier Transform, a product on ${A}$, which is the {\it convolution product} in the classical examples. We prove algebraically Parseval's Identity (Theorem \ref{Parseval}), the Peter-Weyl Theorem (Eq. \ref{Peter}), the Inversion Formula (Eq. \ref{Inversion}), etc. Harmonic Analysis has been developed from an algebraic point of view in finite groups and in the Cartier Duality on multiplicative groups (see \cite{Serre}, \cite{Dm}). A. Van Daele, in the more elaborate algebraic framework of regular multiplier Hopf algebras with invariant functionals, defined the Fourier Transform and proved Plancherel's formula (see \cite{Van}).

In order to compute $w_G$, when $G= Sl_n, Gl_n, O_n, Sp_{2n}$, we consider a system of coordinates in $G$, that is, we
consider $G={\rm Spec}\, A$ as a closed subgroup of a semigroup of
matrices $M_n={\rm Spec}\, B$. Then $A$ is the quotient of $B$ by the
ideal $I$ of the functions of $M_n$ vanishing on $G$. Hence, $A^*$ is a subalgebra of $B^*$ and one has that $k \cdot w_G = A^{*G} = \{ w \in B^{*G} : w(I) = 0 \}$. Moreover, $B^G$ (which is the ring of functions of $M_n/G$), coincides essentially with $B^{*G}$, via the Fourier transform. Finally, we prove that given $w \in B^{*G}$, the
condition $w(I)=0$ is equivalent to $w(I^G)=0$, which is a system of ``homogeneous'' equations, finite in each degree. In the theory of invariants the calculation of the invariant integral  is of great interest, because it yields the calculation of the invariants of any representation, as it is shown.

\section{Preliminary results}

If $G$ is a smooth algebraic group over an algebraically closed field  one may only regard the rational points of $G$ in order to resolve many questions in the theory of linear representations of $G$. In general, for any $k$-group, one may regard the functor of points of the group and its linear representations and in this way $k$-groups and their linear representations may be managed as mere sets, as it is well known.

Let $k$ be a commutative ring with unit. All functors considered
in this paper are functors over the category of commutative
$k$-algebras. Given a $k$-module $E$, the functor ${\bf E}$ defined by ${\bf E}(B) := E \otimes_k B$ is called a quasi-coherent $k$-module. The functors $E \rightsquigarrow {\bf E}$, ${\bf E} \rightsquigarrow {\bf E}(k)$ establish an equivalence between the category of $k$-modules and the category of quasi-coherent $k$-modules (\cite[1.12]{Amel}). In particular, ${\rm Hom}_k ({\bf E},{\bf E'}) = {\rm Hom}_k (E,E')$.

If $F$ and $G$ are functors of $k$-modules, we will denote by ${\bf Hom}_k (F, G)$ the functor of $k$-modules $${\bf Hom}_k (F,G)(B) := {\rm Hom}_B (F_{|B}, G_{|B})$$ where $F_{| B}$ is the functor $F$ restricted to the category of commutative $B$-algebras. We denote $F^* = {\bf Hom}_k (F, {\bf k})$ and given a $k$-module $E$, ${\bf E^*}(B)={\rm Hom}_k(E,B)$. It holds that ${\bf E}^{**} = {\bf E}$  (\cite[1.10]{Amel}).

Let $G={\rm Spec}\, A$ be an affine $k$-group and let $G^{\cdot}$
be the functor of points of $G$, that is, $G^{\cdot}(B) = {\rm
Hom}_{k-alg} (A,B)$.  One has a natural morphism
$G^{\cdot} \to {\bf A}^*$, because $G^{\cdot}(B) = {\rm Hom}_{k-alg} (A,B) \subset {\rm Hom}_k (A,B) = {\bf A}^*(B)$. If $F^*$ is a functor of $k$-modules (resp. $k$-algebras) any morphism of functor of sets (resp. groups) $G^{\cdot}\to F^*$ factorizes uniquely through a morphism of functor of $k$-modules (resp. $k$-algebras) ${\bf A^*}\to F^*$. (\cite[5.3]{Amel}). Therefore, the category of $G$-modules is equal to the category of (quasi-coherent) ${\bf A}^*$-modules (\cite[5.5]{Amel}).

\begin{defi} 
Let $k$ be a field. We will say that an affine $k$-group $G$ is (linearly) semisimple if it is linearly reductive, that is, if
every linear representation of $G$ is completely reducible.\end{defi}

$G$ is semisimple if and only if ${\bf A}^*$ is semisimple, i.e., $A^* = \prod_i A_i^*$, where $A_i^*$ are simple (and finite)
$k$-algebras (\cite[6.8]{Amel}). If $k$ is an algebraically closed field, then $A_i^*$ is an algebra of matrices by Wedderburn's theorem.

\begin{den}
For abbreviation, we sometimes use $g \in G$ or $f \in F$ to denote $g \in G^\cdot(B)$ or $f \in F(B)$ respectively. Given $f \in F(B)$ and a morphism of $k$-algebras $B \to B'$, we still denote by $f$ its image by the morphism $F(B) \to F(B')$
\end{den}

\section{The Fourier Transform}

Assume $k$ is a field.

If $B$ is a $k$-algebra and $E$ is a left $B$-module (resp. a right $B$-module) then $E^*={\rm Hom}_k(E,k)$ is a right $B$-module (resp. a left $B$-module): $( w\cdot b)(e):=w(b\cdot e)$ (resp. $(b\cdot w)(e):=w(e\cdot b)$) for all $w\in E^*$, $e\in E$ and $b\in B$.

Let $B$ be a finite $k$-algebra, let $tr_B\in B^*$ be defined by $tr_B(b):=$ trace of the endomorphism of $B$, $b'\mapsto bb'$, and ley $T_2^B$ be the metric of the trace of $B$, that is, $T_2^B(b,b'):=tr_B(bb')$. Let $\phi_B\colon B\to B^*$ the polarity associated to $T_2^B$. Since $T_2^B(bb',b'')=T_2^B(b',b''b)$ and $T_2^B(b'b,b'')=T_2^B(b',bb'')$, then $\phi_B$ is a morphism of left and right $B$-modules. Obviously, $\phi_B(1)=tr_B$.

Let $B_i$ be finite $k$-algebras and   $A^*:=\prod_i B_i$.  On $\tilde A:=\oplus_i B_i$ one has $tr=\sum_i tr_{B_i}\in \oplus_i B_i^*=A$ and the symmetric metric of the trace $T_2$ (in fact, given $b\in A^*$ and $\tilde b\in \tilde A$ we may define $T_2(\tilde b,b):=tr(\tilde bb)=:T_2(b,\tilde b)$). Again, $$T_2(bb',b'')=T_2(b',b''b), \quad T_2(b'b,b'')=T_2(b',bb'')$$
for all $b\in A^*$, $b',b''\in \tilde A$. The associated polarity $\phi \colon \tilde A=\oplus B_i\to \oplus_i B_i^*=A\subset (\tilde A)^*$ is a morphism of left and right $A^*$-modules. Let $1_i:=(0,\ldots,\overset i1, \ldots, 0)\in B$, then $\phi(1_i)=tr_{B_i}$. Obviously, $\tilde A$ is an $A^*$-module generated by $\left\{1_i\right\}_i$, and $$T_2(1_i,1_j)=tr(1_i\cdot 1_j)=\left\{\begin{array}{ll} \dim_k B_i & \text{ if }  i =j \\ 0 & \text{ otherwise } \end{array} \right.$$

\begin{com} Observe that $\tilde A=\oplus_iB_i=\{b\in A^*=\prod_i B_i\colon \dim_k(A^*\cdot b)<\infty\}$. In functorial words, $\tilde{\bf A}:=\oplus_i {\bf B}_i$ is the maximal quasi-coherent ideal of $\prod_i {\bf B}_i=:\bf A^*$: Let $b = (b_i)_{i} \in \prod_i {\bf B}_i$ and assume $b$ belongs to a quasi-coherent ${\bf A^*}$-submodule of ${\bf A}^*$. The elements $( \ldots, 0, b_i, 0, \ldots) = (\ldots, 0, 1, 0, \ldots) \cdot b$ belong to the quasi-coherent ${\bf A}^*$-submodule $<b>$ generated by $b$. By Proposition \cite[4.7]{Amel}), $<b>$ is a finite $k$-module. Hence, $b_i \neq 0$ only for a finite set of indices $i$, and $b$ belongs to $\tilde{\bf A}$.
\end{com}

\begin{teo} \label{teo}
Let $B_i$ be finite $k$-algebras and   $A^*:=\prod_i B_i$. If
$F\colon A\to A^*$ is a morphism of left and right $A^*$-modules then there exists an element $z\in Z(A^*)$ such that $\phi\circ F=z\cdot$ ($z\cdot (a):=z\cdot a$). Moreover, $F(tr_{B_i})= z\cdot 1_i=z_i$, where $z=(z_i)_i\in \prod_iZ(B_i)=Z(A^*)$.
\end{teo}

\begin{proof} 
Since $F$ is a morphism of $A^*$-modules, $F(A)\subset \tilde A=\oplus B_i$. Hence, $$\phi\circ F\in {\rm Hom}_{A^*\otimes A^*}(A,A)\subset {\rm Hom}_{A^*\otimes A^*}(A^*,A^*)=Z(A^*).$$
Finally, $\phi(F(tr_{B_i}))=z\cdot tr_{B_i}$ and $\phi(z\cdot 1_i)=z\cdot tr_{B_i}$ then $F(tr_{B_i})= z\cdot 1_i=z_i$.
\end{proof}

With the notation of the previous theorem, if $\phi\circ F={\rm Id}$ then
\begin{equation} \label{eq1}  
a(b)=(\phi(F(a)))(b)= T_2(F(a),b)=tr(F(a)b)
\end{equation}
for all $a\in A$ and $b\in A^*$.

Let $G={\rm Spec}\, A$ be a semisimple affine group. One has $A^*
= \prod_{i} A_i^*$, where $A_i^*$ are finite simple algebras. Then we have the metric of the trace on $\tilde A=\oplus_i A_i^*$ and the associated polarity $\phi\colon \tilde A=\oplus_i A_i^*\to \oplus_i A_i=A$.

Let $E$ be a linear representation of a $k$-group $G={\rm Spec}\,
A$. The associated character $\chi_E \in A$ is defined by
$\chi_E(g)=$ trace of the linear endomorphism $E \to E$, $e
\mapsto g \cdot e$, for every $g \in G$ and $e \in E$.

Assume for simplicity that $k$ is an algebraically closed field, then $A_i^*={\rm End}_k (E_i)$. Observe that $tr_{A_i^*}=n_i\cdot \chi_{E_i}$, where $n_i=\dim_k E_i$.

Let  $E_0 = k$ be the trivial representation of $G$ and let $w_G :=
1_0 \in A^*$ be the ``invariant integral on $G$''. The invariant
integral on $G$ is characterized by being $G$-invariant and
normalized, that is, $w_G(1)=1$ (\cite[2.11]{Reynolds}). Given
a functor of $k$-modules $F$ then $H=F^*$ is a functor of $G$-modules if and only if it is a functor of ${\bf A}^*$-modules and the morphism $w_G \cdot: H \to H, h \mapsto w_G \cdot h$ is the unique projection of $G$-modules of $H$ onto $H^G$, in particular $H^G = w_G \cdot H$ (see \cite[2.3,3.3]{Reynolds}).
Given $a\in A$, then $w_G\cdot a\in k = A^G$. Hence, $w_G\cdot a =
(w_G\cdot a)(1)= a(w_G)=w_G(a)$.

One has that $w_G \cdot A_j= 1_0\cdot A_j=0$, if $j\neq 0$  and $w_G \cdot a_0 = a_0$ for all $a_0\in A_0$. Hence, $w_G \cdot \chi_{E_j} = 0$ if $E_j$ is not the trivial representation, and $w_G \cdot \chi_{E_0} = \chi_{E_0} = 1$. Moreover, since $\chi_{E\oplus E'} = \chi_E + \chi_{E'}$, one has
\begin{equation} 
w_G(\chi_E)=w_G \cdot \chi_E = \dim_k E^G .
\end{equation}

Let $* : A \to A$, $a \mapsto a^*$ be the morphism induced  by the
morphism $G \to G$, $g \mapsto g^{-1}$. If $E$ is a linear representation of $G$, we will consider $E^*$ as a left $G$-module by $(g*w)(e)= w(g^{-1}\cdot e)$. One has that $\chi_E^* = \chi_{E^*}$, because the trace of $g^{-1} \in G$ operating on $E$ is equal to the trace of $g$ operating on $E^*$ (which operates by the transposed inverse of $g$).

\begin{teo}\label{Parseval}
Let $G={\rm Spec}\, A$  be a semisimple group and let $w_G \in
A^*$ be its invariant integral. The morphism $$F\colon {A} \to {A}^*,\, \, F(a)= w_G(a^* \cdot -),$$ where $w_G (a^* \cdot -)(a') := w_G(a^* \cdot a')$, is equal to the inverse morphism of $\phi$.
\end{teo}

\begin{proof}
Let us first prove that $F:{\bf A} \to {\bf A}^* ,\, F(a) := w_G(a^* \cdot -)$ is a morphism of left $G$-modules. For every point $g \in G$, $$\aligned w_G((g\cdot a)^*\cdot -)& \overset*=w_G((a^*\cdot g^{-1})\cdot -) \overset{**}=w_G(((a^*\cdot g^{-1})\cdot -)\cdot g)= w_G(a^*\cdot (-\cdot g)) \\ & = g\cdot (w_G(a^*\cdot -))\endaligned$$ where $\overset{*}{=}$ is due to $(g \cdot a)^*(g') = a({g'}^{-1} \cdot g) = a((g^{-1} \cdot g')^{-1}) = (a^* \cdot g^{-1}) (g')$, and $\overset{**}{=}$ is due to $g \cdot w_G = w_G$. Likewise, $F$ is a morphism of right $G$-modules.

Assume $k$ is algebraically closed. Then $A^*=\prod_i {\rm End}_k (E_i)$. By Theorem \ref{teo}, $\phi\circ F = z\cdot$, $z=(z_i)\in \prod_i k=Z(A^*)$ and $F(tr_{{\rm End}_k (E_i)})= z_i\cdot 1_i$.  Hence, $F(\chi_{E_i})=(z_i/n_i)\cdot 1_i$.
Observe that $$F(\chi_{E_i})(\chi_{E_j})  =w_G (\chi_{E_i^* \otimes E_j}) = w_G \cdot \chi_{E_i^* \otimes E_j} = \dim_k {\rm Hom}_G (E_i,E_j) = \delta_{ij}.$$ Hence, $1 = F(\chi_{E_i})(\chi_{E_i})=(z_i/n_i) \cdot 1_i(\chi_{E_i})=(z_i/n_i)\cdot n_i=z_i$ and $\phi\circ F ={\rm Id}$.
\end{proof}

\begin{defi}
Let $G= {\rm Spec}\, A$ be an affine semisimple group and let $w_G \in A^*$ be its invariant integral. The morphism $$ F: A \to A^*, F(a) := w_G(a^* \cdot -),$$ is said to be the Fourier Transform of $G$.
\end{defi}

If we consider on $A$ the metric defined by $T^2(a,a') := w_G(a^*\cdot a')$, then $F$ is the polarity associated to this metric. If we consider the metric $T_2$ on $\tilde{A}=\oplus_i A_i^*$, then $F: A\to \tilde{A}$ is an isometry.
Observe that $F(\chi_{E_i})=(1/n_i)\cdot 1_i$, then $g.c.d.({\rm char}\, k,n_i)=1$.
Observe that

$$T^2(\chi_{E_i},\chi_{E_j})=T_2(1_i/n_i,1_j/n_j)=tr(1_i/n_i\cdot 1_j/n_j)=\left\{\begin{array}{ll} 1 & \text{ if }  i =j \\ 0 & \text{ otherwise } \end{array} \right.$$

 If we denote $w_G = \int_G dg$, then

\begin{equation}\int_G a (g^{-1}) \cdot b(g)\, dg = T^2(a,b) = T_2(F(a), F(b)) = tr(F(a) \cdot F(b)), \end{equation}
(the corresponding classical analytical theorem is known as ``Parseval's identity'', see \cite[2.7.3]{Edwards} and \cite{Gelfand}).

Let $*: A^* \to A^*$ be the transposed morphism of $*: A \to A$. It is easy to check that $*(g)= g^{-1}$ for all $g \in G$. By \cite[2.11]{Reynolds},  $*w_G=w_G$. The Fourier Transform commutes with $*$, that is,  \begin{equation}* \circ F = F \circ *,\end{equation} because given $a,b \in A$ then $((* \circ F)(a))(b) = w_G (a^* \cdot b^*) = w_G (a \cdot b) = ((F \circ *)(a)) (b).$

Given $a \in A$, the diagram $$\xymatrix{  A \ar[r]^-{F} \ar[rd]_-{F(a)} & A^* \ar[d]^-a\\ & k}$$ is commutative, because $(a \circ F)(b) = F(b)(a) = T^2(b,a) = T^2 (a,b) = F(a)(b)$, for all $b \in A$. Let $\pi_1: A^*\to k$ be defined by $\pi_1(w)=w(1)$. Then, $\pi_1\circ F=w_G$.

\begin{com}
Let $\bar{T}^2(a,b):= T^2 (a^*,b) = w_G(a \cdot b)$. It holds:
\begin{enumerate}
\item The polarity associated to $\bar{T}^2$ is $F \circ *$.
\item $\bar{T}^2$ is symmetric.
\item $\bar{T}^2(g \cdot a, g \cdot b) = w_G ((g \cdot a) \cdot (g \cdot b)) = w_G( g \cdot (a \cdot b)) = w_G (a \cdot b) = \bar{T}^2(a,b)$. Likewise, $\bar{T}^2(a,b) = \bar{T}^2(a \cdot g, b \cdot g)$.
\end{enumerate}
\end{com}

\section{Algebraic Harmonic Analysis}

Assume $k$ algebraically closed for simplicity. Let $G= {\rm Spec}\, A$ be an affine semisimple group. Let $\hat G$ be the set of the irreducible linear representations $\rho_i\colon G\to {\rm End}_k (E_i)$ of $G$ (up to isomorphism).

Let $1_i=(0,\ldots,\overset i1, \ldots, 0)\in \prod_{\rho_j \in \hat{G}} {\rm End}_k (E_j) =A^*$. Given $a \in A$ let us calculate the $i$-th-component of $F(a) \in A^* = \prod_i A_i^* = \prod_{\rho_i \in \hat{G}} {\rm End}_k (E_i)$, that is, $1_i \cdot F(a)$. Given $g \in G \subset \bf A^*$, let us denote $(g_{rs}^i) = \rho_i(g)=1_i \cdot g$. Let us see that
\begin{equation} 
1_i \cdot F(a) = \int_G a(g) \cdot (g_{rs}^i)^{-1} dg . \end{equation} 
For all $b \in A$, $(b \cdot 1_i)^*(g) = b((g_{rs}^i)^{-1})$, then $$\aligned \left( \int_G a(g) \cdot (g_{rs}^i)^{-1} dg \right) (b) & = \int_G a(g) \cdot b((g_{rs}^i)^{-1}) dg = \int_G a(g) \cdot (b \cdot 1_i)^*(g) dg \\ & = T^2(b \cdot 1_i, a) = T^2 (b, 1_i \cdot a) \\ & = F(1_i \cdot a)(b) = (1_i \cdot F(a))(b).  \endaligned$$

Given $\rho_i \in \hat{G}$, denote ${\rm End}_k (E_i) = M_{n_i}(k)$ and $w_{rs}^i \in M_{n_i}(k)$ the null matrix except for the $rs$-coefficient which is equal to $1$. Obviously,
$$ T_2(w_{rs}^i, w_{r's'}^{i'}) = \left\{ \begin{array}{cl} n_i & \text{ if } i=i', r=s' \text{ and } s=r' \\ 0 & \text{ otherwise } \end{array} \right.$$
Let $\{ \delta_{rs}^i \}_{i,r,s}$ be the dual basis of $\{ w_{rs}^i \}_{i,r,s}$. Observe that the composite morphism $$ G^\cdot \hookrightarrow {\bf A}^* \stackrel{\delta_{rs}^i}{\to} {\bf k}$$ is equal to $g_{rs}^i$, that is, $\delta_{rs}^i(g) = g_{rs}^i$. $F( \delta_{rs}^i) = w_{sr}^i/n_i$, because $T_2 (F(\delta_{rs}^i), w_{r's'}^{i'}) = \delta_{rs}^i (w_{r's'}^{i'}) = T_2(w_{sr}^i, w_{r's'}^{i'})/n_i$. Hence, 
\begin{equation} \label{Peter} 
T^2(\delta_{rs}^i, \delta_{r's'}^{i'}) =T_2(w_{rs}^i/n_i, w_{r's'}^{i'}/n_{i'})=\left\{ \begin{array}{cl} 1/n_i & \text{ if } i=i', r=s' \text{ and } s=r' \\ 0 & \text{ otherwise } \end{array} \right. 
\end{equation} 
(the corresponding  classical analytical theorem is known as the ``Peter-Weyl Theorem'', see \cite[5.12]{Folland}).

Let $E$ be a $G$-module and let $E_i$ be a simple $G$-module. Let us consider the $G$-module decomposition $E=E'\oplus F$, where $E'$ is the homogeneous component of $E$ isomorphic to $\oplus^n E_i$. Now, we want to compute the morphism $E \to E$ which is the identity over $E'$ and null on $F$. In particular, we could obtain the decomposition of $E$ as a direct sum of homogeneous modules.
Let $1_i=(0,\ldots,\overset i1, \ldots, 0)\in \prod_{\rho_j \in \hat{G}} {\rm End}_k (E_j) =A^*$. We have to calculate the morphism $E \to E,$ $e \mapsto 1_i \cdot e$. Recall that $1_i = F(n_i \cdot \chi_{E_i})$, $n_i = \dim_k E_i$. The dual morphism of the multiplication morphism ${\bf E}^* \otimes {\bf A}^* \to {\bf E}^*$ is the comultiplication morphism $\mu: E \to E \otimes A$. If $\{e'_l\}$ is a basis of $E$, and $\mu(e)= \sum_l e'_l\otimes a_l$, then $g \cdot e = \sum_l a_l(g)e'_l$ for all $g \in G$. Hence,
\begin{equation}
\aligned  1_i \cdot e & = \sum_l a_l(1_i) \cdot e'_l = \sum_l a_l(F(n_i \cdot \chi_{E_i})) \cdot e'_l = \sum_l n_i \cdot w_G(a_l \cdot \chi_{E^*_i}) \cdot e'_l \\ & = n_i \cdot \sum_l ( \int_G a_l(g) \cdot \chi_{E^*_i}(g) \, dg ) \cdot e'_l.\endaligned
\end{equation}

Now let us prove a formula from the classical harmonic analysis with strong applications in the solution of differential equations (see \cite[ch. 2, 2.1.(3)]{Gelf-Ch})

\begin{pro}
Let $G= {\rm Spec}\, A$ be a semisimple group, let $D$ be a left $G$-invariant vector field on $G$ and let $D_e$ be the value of the vector field at the identity element $e \in G$. Then  \begin{equation} 
F(D(a)) = D_e \cdot F(a) \quad \forall a \in A. 
\end{equation}
\end{pro}

\begin{proof}
It holds that $D(a) = D_e \cdot a$, because $$ (D_e \cdot a)(g) = a ( g \cdot D_e) = a (D_g) = D_g(a) = D(a)(g)$$ for all $g \in G$. Therefore, $F(D(a)) = F(D_e \cdot a )= D_e \cdot F(a)$.
\end{proof}

\begin{den}
Given an affine scheme $X={\rm Spec}\, A$, we will denote $A_X=A$.
\end{den}

\begin{pro} 
Let $G={\rm Spec}\, A$ be a semisimple group, $H\overset i\subset G$ a normal subgroup and $G \overset{\pi}{\to} G/H$ the quotient morphism. Let $i^* : A_{G} \to A_H$ and $\pi: A_G^* \to A_{G/H}^*$ be the natural morphisms. Then, with the obvious notation,
\begin{equation}
w_H(i^*(a))=tr_{G/H}(\pi(F_G(a))), 
\end{equation} 
for all $a \in A_G$. (In the classical harmonic analysis this expression can be understood as the Poisson summation formula, see \cite{Gelfand}).
\end{pro}

\begin{proof} 
The set of the irreducible representations of $G/H$ is equal to the set of the irreducible representations of $G$ which are $H$-invariant. The natural projection $A^*_G=\prod_i {\rm End}_k (E_i) \to \prod_{E_i=E_i^H} {\rm End}_k(E_i)=A_{G/H}^*$ coincides with $\pi$.

The diagram
$$\xymatrix{ A_G =\oplus_i {\rm End}_k (E_i)^* \ar[r]^-{F_G} \ar[d]^-{w_H\cdot} &
\prod_i {\rm End}_k (E_i) \ar[d]^-{w_H\cdot} \ar@{=}[r] & A_G^*  \ar[d]^-{w_H\cdot}
\ar@{=}[r] & A^*_G \ar[d]^-{\pi}
\\
A_G^H=\oplus_{E_i=E_i^H} {\rm End}_k (E_i)^*\ar[r]^-{F_{G/H}} &
\prod_{E_i=E_i^H} {\rm
End}_k(E_i) \ar@{=}[r]& A_{G}^{*H}\ar@{=}[r] & A_{G/H}^*} $$ 
is commutative. Then
$$tr_{G/H}(\pi(F_G(a)))=tr_{G/H}(F_{G/H}(w_H\cdot
a))=(w_{H}\cdot a)(1)=a(w_H)=w_H(i^*(a)) .$$
\end{proof}

The image of the morphism $F: A \hookrightarrow A^*$ is a bilateral ideal. Then it is a subring, although without unit if $\dim_k A=\infty$, because $(\ldots, 1, 1, 1, \ldots) \not\in \oplus A_i^* = {\rm Im}\, F$.

\begin{defi}\label{convolution}
The product of the subring ${\rm Im}\,F$ induces a product on $A$, through the identification $A \stackrel{F}{\backsimeq} {\rm Im}\,F$. This product is called the convolution product.
\end{defi}

Let $a, b \in A$, and let us denote by $*$ the convolution product. Then for all $x \in G$ 
\begin{equation}
\aligned(a*b)(x)  & = T_2(F(a*b), x) = T_2(F(a) \cdot F(b), x) = T_2 (F(a), F(b) \cdot x) \\ & = T_2(F(a), F(b \cdot x) ) = T^2 (a, b \cdot x) = \int_G a(g^{-1}) \cdot b(x \cdot g) dg  .\endaligned
\end{equation} 
In particular, given the identity element $1 \in G$, $$ (a * b)(1)= \int_G a(g^{-1}) \cdot b(g) dg = T^2(a,b).$$

\vskip1pt {\bf Inversion Formula} \vskip1pt

Let $G = {\rm Spec}\, A$ be an affine semisimple group. Let $\tilde A={\rm Im}\,F = \oplus_i A_i^*$..It is easy to check that the set of maximal bilateral ideals of $\tilde{A}$ is equal to $\hat{G}$.

\begin{defi}
We will say that $\tilde{A} \subseteq A^*$ is the algebra of functions of $\hat{G}$.
\end{defi}

Let $e \in G$ be the identity element of $G$, that is, the unit of $A^*$. The diagram $$\xymatrix{  A \ar[r]^-{F} \ar[rd]_-{e} & \tilde{A} \ar[d]^-{tr} \\ & k}$$ is commutative, because $tr(F(a)) = tr (F(a) \cdot e) = T_2 ( F(a), e) = a(e)$, that is, $tr \circ F = e$. Moreover, $tr \circ ( F \circ *) = e \circ * = e$ and  \begin{equation}tr \circ * = e \circ F^{-1} \circ * = (e \circ * ) \circ F = e \circ F^{-1} = tr. \end{equation}

We know that ${\rm Hom}_G(A,k)={\rm Hom}_{A^*}(A,k)=k\cdot w_G$. We want to prove that $tr$ is the ``invariant integral'' on $\hat G$.

$A$ is a left (and right) $A$-module, then $A^*$ is an $A$-module, whose product we denote by $ \bar{\cdot}$. Observe that $F(a) = w_G(a^* \cdot -) = a^* \bar{\cdot} w_G$, then $\tilde{A} = A \bar{\cdot} w_G$ and $\tilde{A}$ is an $A$-submodule of $A^*$. The isomorphism $F \circ * : A \to \tilde{A}$ is a morphism of $A$-modules because $( F \circ * )(a) = a \bar{\cdot} w_G$.

\begin{pro}
Let $k$ be an $A$-module as follows: $a \bar{\cdot} \lambda := a(e) \cdot \lambda $, for all $a \in A$ and $\lambda \in k$. It holds $${\rm Hom}_A(\tilde{A}, k) = k \cdot tr .$$
Hence, $tr: \tilde{A} \to k$ is the unique morphism of $A$-modules such that $tr(w_G)=1$.
\end{pro}

\begin{proof}
Obviously, ${\rm Hom}_A(A, k)= k \cdot e$ and $F \circ * :A \to \tilde{A}$ is an isomorphism of $A$-modules. Hence, ${\rm Hom}_A (\tilde{A}, k) = k \cdot e \circ (F \circ *)^{-1} = k \cdot tr$.

Finally, $tr(w_G) = tr( F(1)) = e(1)= 1$.
\end{proof}

If $G= {\rm Spec}\,A$ is a finite (and etal\'e) commutative group and $g.c.d.({\rm car}\, k, \#G)=1$, then $\hat{G}={\rm Spec}\, A^*$. Moreover, $\tilde{A} = A^*$ and $tr = \# G \cdot w_{\hat{G}}$.

By abuse of notation, we write $tr = \int_{\hat{G}}$.

\begin{defi}
We will call the morphism $$F_{\hat{G}} : \tilde{A} \to A \subset (\tilde{A})^*, \; F_{\hat{G}} (w) := \int_{\hat{G}} (w^* \cdot -)$$ where $F_{\hat{G}} (w)(w')= \int_{\hat{G}} (w^* \cdot w') = tr (w^* \cdot w')$ for all $w \in \tilde{A}$, $w' \in A^*$, the Fourier Transform of $\hat{G}$.
\end{defi}

\begin{teo}
Let $G = {\rm Spec}\, A$ be an affine semisimple $k$-group. The diagram $$\xymatrix{ A \ar[dr]_-{F_G} \ar[rr]^-{*} & &  A \\ & \tilde{A} \ar[ur]_-{F_{\hat{G}}} & }$$ is commutative.
\end{teo}

\begin{proof}
It holds $(F_{\hat{G}}(F_G(a)))(w) = tr (F_G(a)^* \cdot w) = tr(F_G(a^*) \cdot w ) = a^*(w)$, by Equation (\ref{eq1}), for all $a \in A$  and $w \in A^*$.
\end{proof}

Given $w = (w_i)_i \in \oplus_{\rho_i \in \hat{G}} {\rm End}_k(E_i) = \tilde{A}$, we denote $\int_{\hat{G}} w = \int_{\hat{G}} w_i d \rho_i$. Then we can rewrite the previous theorem as follows  \begin{equation} \label{Inversion} 
a(g) =a^*(g^*) = (F_{\hat{G}}(F_G(a)))(g^*) = \int_{\hat{G}} F_G(a)^* \cdot g^* = \int_{\hat{G}} F_G(a)_i \cdot (g^i_{rs}) d\rho_i, 
\end{equation} 
which allows to recover a function on $G$ by means of its Fourier Transform (the corresponding classical analytical expression is known as the ``Inversion Formula'', see \cite[2.8.8]{Edwards}, \cite[5.15]{Folland} and \cite{Gelfand}).

\section{Invariant integral on $Sl_n$, $Gl_n$, $O_n$ and $Sp_{2n}$}

Let $k$ be a field of characteristic zero. The groups $Gl_n$,
$Sl_n$, $O_n$ and $Sp_{2n}$ are semisimples, so they have an invariant integral. This section is devoted to the explicit calculation of the invariant integral on the groups $Gl_n$, $Sl_n$, $O_n$ and $Sp_{2n}$.

Let us consider the affine algebraic $k$-variety $M_n = {\rm
End}_k (E)$, whose points with values in a $k$-algebra $B$ is the semigroup of square matrices of order $n$ with coefficients in $B$. Its ring of functions is $A_{M_n} = \underset{n \in \N}{\oplus} S^n
({\rm End}_k (E)^*)$. One has  ${ A}_{M_n}^* = {\underset{n \in \N}{\prod}} ( {\rm End}_k ({ E}) \otimes_k \stackrel{n}{\ldots} \otimes_k {\rm End}_k ({E}) )^{S_n}$ and the morphism $$ M_n^{\cdot} \to {\bf A}_{M_n}^*, \; \tau \mapsto (1, \tau, \tau \otimes \tau, \tau \otimes \tau \otimes \tau, \ldots ) . $$
There exists a unique structure of functor algebras on ${\bf A}_{M_n}^*$ such that this morphism is a morphism of functor of semigroups. Specifically, $A^*_{M_n}$ is the direct product of the algebras $({\rm End}_k(E)^{\otimes n})^{S_n}\subset
{\rm End}_k(E)^{\otimes n}$ and given $T_1\otimes\cdots \otimes T_n,S_1\otimes\cdots \otimes S_n\in {\rm End}_k(E)^{\otimes n}$ then $(T_1\otimes\cdots \otimes T_n)\cdot (S_1\otimes\cdots \otimes S_n):=
T_1S_1\otimes\cdots \otimes T_nS_n$.

The natural action of ${\rm End}_k (E)^\cdot$ on $E \otimes
\overset{m}{\ldots} \otimes E$ extends to a unique structure of
${A}_{M_n}^*$-module. It consists of the projection of $\underset{m}{\prod} S^m {\rm End}_k ({ E})$ onto the $m$-th factor, $S^m {\rm End}_k ({E})$, and the action of $S^m {\rm End}_k ({E})$ on $E \otimes \overset{m}{\ldots} \otimes E$ via its inclusion in ${\rm End}_k (E) \otimes \overset{m}{\ldots} \otimes {\rm End}_k (E)$, that is, $$(g_1 \cdot \ldots \cdot g_m) \cdot (v_1 \otimes \ldots \otimes v_m) = \frac{1}{m!} \cdot \sum_{\sigma \in S_m} g_{\sigma(1)}(v_1) \otimes \ldots \otimes g_{\sigma(m)} (v_n).$$

The natural isomorphism ${\rm End}_k (E)^* \to {\rm End}_k (E)$ induces, by taking symmetric algebras, a morphism $\varphi : A_{M_n} \hookrightarrow A^*_{M_n}$ of left and right $A_{M_n}^*$-modules. The morphism $\varphi$ coincides with the Fourier Transform $F$, up to an invertible factor of the centre $Z(A^*_{M_n})$, by Theorem \ref{teo}.

\vskip3pt{\bf Invariant integral on $Sl_n$.} \vskip3pt

Let $A_{Sl_n}= k[x_{11}, \ldots, x_{nn}]/(det (x_{ij}) -1)$ be the
ring of functions of the special linear group.  Recall that $ ({ A}_{Sl_n}^*)^{Sl_n} = { k} \cdot w_{Sl_n}$. Let us compute the invariants of ${A}_{Sl_n}^*$ by $Sl_n$.

From the inclusion $Sl_n \subset M_n $ one obtains the injective
morphism $A_{Sl_n}^* \subset A_{M_n}^*$. We will first
compute the invariants of $A_{M_n}^*$ by the action of $Sl_n$ and then we will compute the ones belonging to $A_{Sl_n}^*$. Since $A_{M_n}$ is a semisimple $Sl_n$-module, it splits into a direct sum of $Sl_n$-modules $A_{M_n}= A_{M_n}^{Sl_n} \oplus B$. Obviously  $(A^*_{M_n})^{Sl_n}=(A_{M_n}^{Sl_n})^*$ and $A_{M_n}^* = (A_{M_n}^*)^{Sl_n} \times B^*$. The morphism $$\varphi: {\bf A}_{M_n} = {\bf A}_{M_n}^{Sl_n} \oplus {\bf B} \to ({\bf A}_{M_n}^*)^{Sl_n} \times {\bf B}^*={\bf A}_{M_n}^*$$ is a
morphism of $Sl_n$-modules, so that $\varphi ( A_{M_n}^{Sl_n} ) \subseteq ( A_{M_n}^*)^{Sl_n}$ and $\varphi (B) \subseteq B^*$.

\begin{den} 
Given a $k$-vector space $E$ we will regard the topology in $E^*$ whose closed sets are $\{(E/V)^*:=\{w\in E^*\colon w(V)=0\}\}_{V\subseteq E}$.
\end{den}

Since the closure $ \overline{\varphi ( A_{M_n} )}$ of $\varphi ( A_{M_n} )$ is $A_{M_n}^*$, one has that $\overline{ \varphi ( A_{M_n}^{Sl_n} ) } = (A_{M_n}^*)^{Sl_n}$. Then, let us compute $A_{M_n}^{Sl_n}$. From the exact sequence of groups $$\begin{matrix} 1 & \to & Sl_n & \subset & Gl_n & \to & Gl_n/Sl_n = G_m & \to & 1 \\ & & & &T & \mapsto & det(T) & & \end{matrix}$$ it follows easily that ${ k[x_{11}, \ldots, x_{nn}]}^{Sl_n} = { k[det(x_{ij})]}$ (which is well known, see \cite[II.0.9]{Popov}) and,
therefore, $${( k[x_{11}, \ldots, x_{nn}]^*)}^{Sl_n} = k
\times k \cdot \varphi(det(x_{ij})) \times \ldots, \times k \cdot \varphi(det(x_{ij})^r) \times \ldots $$
Let us denote by $\delta_{ij}$ the matrix of null coefficients,
except for the $ij$-th coefficient that is $1$. Observe that $\delta_{ij} \in A^*_{M_n(k)}$ coincides with $\dfrac{\partial}{\partial x_{ij}}\left. \right|_0$. One has that $$\varphi(x_{ij})= \dfrac{\partial}{\partial x_{ji}} \left. \right|_0 ,\quad \varphi(x_{i_1 j_1} \cdot \ldots \cdot x_{i_m j_m}) = \dfrac{1}{ m!} \dfrac{\partial}{\partial x_{j_1 i_1}} \cdot \ldots \cdot \dfrac{\partial}{\partial x_{j_m i_m}} \left.\right|_0 . $$ Since $det(x_{ij})^r$ is an homogeneous polynomial of $rn$-th degree, $$ \varphi(det(x_{ij})^r) = \dfrac{1}{(rn)!} det^r \left( \dfrac{\partial}{\partial x_{ij}} \right) \left. \right|_0 . $$ Let $D= det \left( \dfrac{\partial}{\partial x_{ij}} \right)$ be the Cayley operator, and let us denote $D_0^r=D^r \left. \right|_0$. Then
$${(k[x_{11}, \ldots, x_{nn}]^*)}^{Sl_n} = k \times k \cdot D_0 \times \ldots \times k \cdot D_0^r \times \ldots$$

Let us compute now the  $\tilde w \in {(k[x_{11}, \ldots, x_{nn}])^*}^{Sl_n}$ vanishing on the ideal $I=(det(x_{ij}) - 1)$. Since $\tilde w$ is $Sl_n$-invariant, one has that $\tilde{w} \cdot w_{Sl_n} = \tilde{w}$. Therefore, $\tilde{w}(I)=0$ if and only if $\tilde{w} (w_{Sl_n} \cdot I) = \tilde{w} (I^{Sl_n}) = 0$. $I^{Sl_n} = \langle det^n(x_{ij}) - det^{n-1}(x_{ij}) \rangle_n \subset k[det(x_{ij})]$. Hence, if $\tilde{w}$ vanishes on the functions $$det^n(x_{ij}) - det^{n-1}(x_{ij}),$$ for every $n\geq 1$, and $\tilde{w}(1) = 1$, then $\tilde{w} = w_{Sl_n}$. Consequently,
 \begin{equation} { w_{Sl_n} = \underset{i}{\sum} \dfrac{D_0^i}{D_0^i
(det^i(x_{ij}))} } . \end{equation} It only remains to determine the value of $D_0^i (det^i(x_{ij})) \in k$.

\begin{lema}{\rm (}\cite[2.1]{Dolgachev}{\rm )}
$D(det^r(x_{ij})) = \mu_{r} det^{r-1}(x_{ij})$, where $\mu_{r} = r
\cdot (r+1) \cdot \ldots \cdot (r+n-1) =
\dfrac{(r+n-1)!}{(r-1)!}$.
\end{lema}

Then  
\begin{equation}
D^r (det^r(x_{ij}))= \mu_r \cdot \mu_{r-1} \cdot \ldots \cdot \mu_1 .
\end{equation}

Every linear representation of $Sl_n$ is a submodule of a direct
sum of $A_{Sl_n}$ (the regular representation). Moreover, $A_{Sl_n}$ is a quotient of the ring of functions of $M_n$. Finally, the ring of
functions of $M_n = {\rm End}_k (E)$ is included in a direct sum
of $E \otimes \overset{m}{\ldots} \otimes E$. Let us compute, then, the invariants of $Sl_n$ acting on these vector spaces.

\begin{pro}{\rm (}\cite[Th. 19.2]{carlos} {\rm )}
Let $Sl_n$ be the special linear group of an $n$-dimensional vector space $E$. Let us consider the natural action of $Sl_n$
on $E \otimes \overset{m}{\ldots} \otimes E$, $g \cdot (v_1
\otimes \ldots \otimes v_m) = g \cdot v_1 \otimes \ldots \otimes g
\cdot v_m$. It holds that:
\begin{enumerate}
\item $(E \otimes \overset{n}{\ldots} \otimes E)^{Sl_n} =
\Lambda^n E$.

\item $(E \otimes \overset{nm}{\ldots} \otimes E)^{Sl_n} =
\sum_{\sigma \in S_{nm}} \sigma ( \Lambda^n E \otimes
\overset{m}{\ldots} \otimes \Lambda^n E)$, where $\sigma \in
S_{nm}$ acts on $E \otimes \overset{nm}{\ldots} \otimes E$
by permuting the factors.

\item $(E \otimes \overset{m}{\ldots} \otimes E)^{Sl_n}=0$ if $m$
is not a multiple of $n$.
\end{enumerate}
\end{pro}

\begin{proof} $ $
\begin{enumerate}
\item We must calculate $w_{Sl_n} \cdot (E \otimes
\overset{n}{\ldots} \otimes E) = D_0 \cdot (E \otimes
\overset{n}{\ldots} \otimes E)$. Fixed a basis $\{ e_1,
\ldots, e_n \}$ of $E$, let us observe that $\frac{\partial}{\partial
x_{ij}}_{|0}$ corresponds to the matrix $\delta_{ij}$ that maps
the vector $e_j$ to $e_i$ and the rest of the $e_k$ to zero. Then
it is clear that $D_0 \cdot e_{i_1} \otimes \ldots \otimes e_{i_n}
= \frac{1}{n!} e_{i_1} \wedge \ldots \wedge e_{i_n}$ and $(1)$ is proved.

\item The $r=nm$-th component of $w_{Sl_n}$ is, up to
scalar, $D^m$; that is, it coincides with
$\underset{\sigma\in S_{nm}}{\sum} \sigma \circ (D \otimes
\overset{m}{\ldots} \otimes D) \circ \sigma^{-1}$, up to scalar. Then, $$(E \otimes \overset{nm}{\ldots} \otimes E)^{Sl_n} \subseteq
\sum_{\sigma\in S_{nm}} \sigma(\Lambda^n E \otimes
\overset{m}{\ldots} \otimes \Lambda^n E) .$$ The inverse inclusion
is obvious.

\item $w_{Sl_n} \in \underset{r}{\prod} S^r {\rm End}_k (E)$ and
its $r$-th ``components'' are null when $r$ is not a multiple of $n$.
\end{enumerate}
\end{proof}

We can compute the dimension of $(E \otimes \overset{nm}{\ldots}
\otimes E)^{Sl_n}$:

$$\aligned \dim_k & (E \otimes \overset{nm}{\ldots} \otimes E)^{Sl_n}  = w_{Sl_n} \cdot \chi_{E \otimes \overset{nm}{\ldots} \otimes E} = \frac{D^m}{D^m(\det^m)} (\chi_E^{nm}) \\ & =
\frac{D^m}{D^m(\det^m)} ((x_{11}+ \overset{n}{\ldots} +
x_{nn})^{nm}) = \frac{D^m}{D^m(\det^m)} (x_{11}^m \cdot \ldots
\cdot x_{nn}^m \cdot \frac{(mn)!}{m!^n}) \\ & =
\frac{(mn)!}{D^m(\det^m)} .\endaligned$$

\vskip2pt {\bf Invariant integral on $Gl_n$.} \vskip3pt

Let $A_{Gl_n}= k[x_{11}, \ldots, x_{nn}, \frac{1}{det(x_{ij})}]$
be the ring of functions of the linear group. One has that $(A_{Gl_n}^*)^{Gl_n} = k \cdot w_{Gl_n}$ and $(A_{Gl_n}^*)^{Gl_n} = ((A_{Gl_n}^*)^{G_m})^{Sl_n}$, where $G_m$ is the multiplication group. First we will compute $(A_{Gl_n}^*)^{G_m} = (A_{Gl_n}^{G_m})^*$ and then we will look for the $Sl_n$-invariant ones among them.

$A_{Gl_n}$ is a $\mathbb{Z}$-graded algebra, whose $i$-th component we denote by $A_i$. Given $\lambda \in G_m$ and $a_i \in A_i$, $\lambda * a_i = \lambda^i \cdot a_i$. Therefore, $A_{Gl_n}^{G_m} = A_0$ and $w \in (A_{Gl_n}^*)^{G_m} $ if and only if it factors via the obvious quotient $A_{Gl_n} \to A_0$. This quotient morphism is a morphism of $Gl_n$-modules. Now we must compute the linear forms $w : A_0 \to k$ that are $Sl_n$-invariant. One has that

$$ A_0 = \underset{r \in \N}{\bigcup} \dfrac{A^r}{det^r(x_{ij})},\qquad A^r:= \{ \mbox{ homogeneous polynomials of $n \cdot r$-th degree } \}.$$ The morphism $w_r = w \circ det^{-r}(x_{ij}) \cdot : A^r \to k$ is a morphism of $Sl_n$-modules, that is, is $Sl_n$-invariant. Since $w_r \in (A^r)^* = (S^{r \cdot n} {\rm End}_k (E)^*)^* = S^{r \cdot n} {\rm End}_k (E) \subset \underset{i}{\prod} S^i {\rm End}_k (E) = k[x_{11}, \ldots, x_{nn}]^*$, and it is $Sl_n$-invariant, then it must be $w_r = \alpha_r \cdot D^r$. If we ask for $w(1)=1$, then it must be $1 = w_r (det^r(x_{ij})) = \alpha_r \cdot D^r (det^r(x_{ij}))$, because $1 = \dfrac{ det^r(x_{ij})}{det^r(x_{ij})} \in A_{Gl_n}^{G_m}$. Consequently, $\alpha_r = \dfrac{1}{D^r (det^r(x_{ij}))}$ and $$w_r = \dfrac{D^r}{ D^r (det^r(x_{ij}))} = \dfrac{D^r}{\mu_r \cdot \ldots \cdot \mu_1} .$$ In conclusion, we have determined\footnote{Marcel B\"{o}kstedt checks in \textit{``Notes on Geometric Invariant Theory''} (available at \texttt{http://home.imf.au.dk/marcel/GIT/GIT.ps}) that the
integral thus defined is the Reynolds operator of the linear
group, and he states that Cayley, in some sense, had already checked
it.} the invariant integral $w_{Gl_n}$ on $Gl_n$ as a linear form over $A_{Gl_n}$:
\begin{equation}
\aligned w_{Gl_n} \left( \dfrac{p(x_{ij})}{det^s(x_{ij})}
\right) & = w_{Gl_n} \left( \dfrac{\ldots + p_{n \cdot s}(x_{ij})+
\ldots }{det^s(x_{ij})} \right) = w_s ( p_{n \cdot s}(x_{ij})) \\
& = \dfrac{D^s (p_{n \cdot s}(x_{ij}))} {D^s (det^s(x_{ij}))} .
\endaligned 
\end{equation}

\vskip1pt {\bf Invariant integral on $O_n$.} \vskip3pt

Let $T_2$ be a non-singular symmetric metric on a vector space
$E$ of dimension $n$. Let $O_n$ be the subgroup of the linear
group of the symmetries of $T_2$. In the algebraic variety $S^2
{\rm E}^*$ of the symmetric metrics, regardless of the basis of
$E$ chosen, we can define (up to a constant multiplicative factor)
the function $det$ that assigns to each metrics its determinant.
So, we can consider the open set $S^2 {\rm E}^* - (det)_0$. The
sequence of morphisms of varieties

$$\begin{matrix} 1 & \to & O_n & \to & Gl(E) & \to & S^2 {\rm E}^* -
(det)_0 & \to & 1 \\ & & & & & & & & \\ & & & & S & \mapsto & S^t
\circ T_2 \circ S & & \end{matrix}$$ shows that $S^2 {\rm
E}^*-(det)_0$ is the quotient variety of $Gl(E)$ by the orthogonal
subgroup $O_n$ (acting $O_n$ on the left on $Gl(E)$). Fixing a
basis in $E$, we will say that $k[x_{11},\ldots,x_{nn},
\dfrac{1}{det(x_{ij})}]$ is the ring of functions of $Gl(E)$ and
$k[y_{i\leq j}, \dfrac{1}{det(y_{ij})}]$ is the ring of functions
of $S^2 {\rm E}^*-(det)_0$. One has the induced morphism of rings $$ \begin{matrix} k[y_{i\leq j}, \dfrac{1}{det(y_{ij})}]
& \hookrightarrow & k[x_{11},\ldots,x_{nn},
\dfrac{1}{det(x_{ij})}] \\ & & \\ y_{rs} & \mapsto & [ (x_{ij} )^t
\circ T_2 \circ (x_{ij}) ]_{rs} \\ & & \\ det(y_{ij}) & \mapsto &
det(x_{ij})^2 \cdot det\, T_2 . \end{matrix}$$ The functions of $Gl(E)$ invariant by $O_n$ identify with the functions of $S^2 {\rm E}^*-(det)_0$. Therefore, via the morphism of varieties ${\rm End}_k (E) \to S^2 {\rm E}^*$, $S\mapsto S^t \circ T_2 \circ S$, the functions of $S^2 {\rm E}^*$ identify with the functions of ${\rm End}_k (E)$ that are (right) invariant by $O_n$.

Let us express these equations without fixing basis. We have
defined the morphism

$$\begin{array}{ccl} {\rm End}_k (E) = {\rm E}^* \otimes {\rm E} &
\longrightarrow & S^2 {\rm E}^* \\ & & \\ w \otimes e & \mapsto &
C^{1,2}_{2,3} ( w \otimes e \otimes T_2 \otimes e \otimes w ) \\ &
& = T_2(e,e) \cdot w \otimes w \end{array} $$ that induces a
morphism between the rings of functions $S^{\cdot} (S^2 E) \to
S^{\cdot} ({\rm End}_k (E)^*)$, that is expressed explicitly as
follows: $$ \begin{array}{ccl} S^m (S^2 E) & \longrightarrow & S^{2m} ({\rm End}_k (E)^*) \\ & & \\ s_1 \cdot \ldots \cdot s_m & \mapsto & \overline{ T_2 \otimes \overset{m}{\ldots} \otimes T_2 \otimes s_1 \otimes \ldots \otimes s_m} \end{array}$$ (we think of $S^{2m} ({\rm End}_k (E)^*)$ as a quotient of $(E^* \otimes E) \otimes \overset{2m}{\ldots} \otimes (E^* \otimes E) = E^* \otimes
\overset{2m}{\ldots} \otimes E^* \otimes E \otimes
\overset{2m}{\ldots} \otimes E$). Equivalently, the left $O_n$-invariant functions of the variety ${\rm End}_k ({\rm E})$ are the direct sum of the images of the morphisms

$$\begin{array}{ccl} S^m (S^2 E^*) & \longrightarrow &  S^{2m}
({\rm End}_k (E)^*) \\ & & \\ \omega_1 \cdot \ldots \cdot \omega_m
& \mapsto & \overline{\omega_1 \otimes \ldots \otimes \omega_m
\otimes  T^2 \otimes \overset{m}{\ldots} \otimes T^2}
\end{array}$$ (we think of $S^{2m}({\rm End}_k (E)^*)$ as a quotient
of $E^{*2m} \otimes E^{2m}$). Therefore, the invariants of
$S^{2m}({\rm End}_k (E)^*)$ by the left and right action of $O_n$ are

$$\aligned  & \langle \overline{\sigma( T_2 \otimes \overset{m}{\ldots}
\otimes T_2) \otimes \sigma' ( T^2 \otimes \overset{m}{\ldots}
\otimes T^2)} \rangle_{\sigma,\sigma' \in S_{2m}} \\ & = \langle
\overline{(T_2 \otimes \overset{m}{\ldots} \otimes T_2) \otimes
\sigma'(T^2 \otimes \overset{m}{\ldots} \otimes T^2)}
\rangle_{\sigma'\in S_{2m}} . \endaligned$$

Let $A_{O_n}$ be the ring of functions of $O_n$ and let $w_{O_n} \in
A^*_{O_n} \subset A^*_{M_n} = \underset{r}{\prod} S^r {\rm End}_k (E)$ be the invariant integral on $O_n$. The $r$-th component $[w_{O_n}]_r$ of $w_{O_n}$ is
$$\begin{array}{lll}
\,[w_{O_n}]_r = & \underset{\sigma \in S_{2m}}{\sum} \lambda_{\sigma} \cdot
\overline{(T_2 \otimes \overset m \cdots \otimes T_2) \otimes
\sigma( T^2 \otimes \overset{m}{\ldots} \otimes T^2)}, &
\text{ if } r=2m, \\ \,[ w_{O_n} ]_r
= & 0, & \text{ if } r=2m+1. \end{array}
$$

\begin{pro}{\rm (}\cite[Th. 4.3.3]{Goodman}{\rm )}
Let us consider the natural action of $O_n$ on $E \otimes
\overset{r}{\ldots} \otimes E$, $g \cdot (e_1 \otimes \cdots
\otimes e_r) = g \cdot e_1 \otimes \ldots \otimes g \cdot e_r$. It
holds that:
\begin{enumerate}
\item $ (E \otimes \overset{2m+1}{\ldots} \otimes E)^{O_n} = 0$.
\item $ (E \otimes \overset{2m}{\ldots} \otimes E)^{O_n} = \langle
\sigma(T^2 \otimes \overset{m}{\ldots} \otimes T^2)
\rangle_{\sigma \in S_{2m}}$.
\end{enumerate}
\end{pro}

\begin{proof} $ $
\begin{enumerate}
\item $w_{O_n} \cdot (E \otimes \overset{2m+1}{\ldots} \otimes E)
= [w_{O_n}]_{2m+1} \cdot E^{\otimes 2m+1} = 0 \cdot E^{\otimes
2m+1}=0$.

\item One has $[w_{O_n}]_{2m} = \underset{\sigma\in S_{2m}}{\sum}
\lambda_{\sigma} \cdot \overline{(T_2 \otimes \overset{m}{\ldots}
\otimes T_2) \otimes \sigma(T^2 \otimes \overset{m}{\ldots}
\otimes T^2)}$ and

$$\aligned & \overline{(T_2 \otimes \overset{m}{\ldots} \otimes T_2) \otimes \sigma(T^2\otimes \overset{m}{\ldots} \otimes T^2)} \\
& =\frac{1}{(2m)!} \sum_{\sigma'\in S_{2m}} {\sigma'(T_2 \otimes
\overset{m}{\ldots} \otimes T_2) \otimes \sigma' \sigma (T^2
\otimes \overset{m}{\ldots} \otimes T^2)} \endaligned$$ via the
inclusion $S^{2m} {\rm End}_k (E) \subset E^{*\otimes 2m} \otimes
E^{\otimes 2m}$. Moreover, $E^{*\otimes 2m} \otimes E^{\otimes
2m}$ acts on $E^{\otimes 2m}$ by contracting each linear form
with the corresponding vector. Therefore, 
\begin{equation*} 
\quad \qquad (E^{\otimes 2m})^{O_n} = w_{O_n} \cdot E^{\otimes 2m} = [w_{O_n}]_{2m} \cdot E^{\otimes 2m} \subseteq \langle \sigma(T^2 \otimes \overset{m}{\ldots} \otimes T^2) \rangle_{\sigma\in S_{2m}} . 
\end{equation*} 
The inverse inclusion is obvious.
\end{enumerate}
\end{proof}

Let us consider the morphism $S^2 E \hookrightarrow {\rm End}_k
(E)$, ${T'}^2 \mapsto {T'}^2 \circ T_2$, that assigns to every
metric ${T'}^2$ the endomorphism associated to the pair of metrics
$T^2,{T'}^2$. Two metrics are isometrics (with regard to
$T_2$) if and only if their associated endomorphisms are equivalent, and every endomorphism (up to conjugation) is the endomorphism associated to a symmetric metric and $T_2$ (\cite{Ermolaev}).
As a result one has that the invariant functions of ${\rm
End}_k(E)$ (by the action by conjugation of the linear group) are invariant functions of $S^2 E$ by the orthogonal group. Conversely, let us see that $A^{O_n}_{S^2 E} \subseteq A^{Gl_n}_{{\rm
End}_k(E)}$. Let $d(\lambda_{ij})$ be the discriminant of the
characteristic polynomial of the matrix $(\lambda_{ij})$ and let
$U = {\rm End}_k (E) - (d)_0$ be the open subset of ${\rm End}_k
(E)$ of the diagonalizable endomorphisms with different
eigenvalues. It is clear that $( U \cap S^2 E ) / O_n = U / Gl_n =
{\rm Spec}\, k[a_1, \ldots, a_n]_d$, where $a_s(\lambda_{ij})$ are
the coefficients of the characteristic polynomial of the matrix
$(\lambda_{ij})$. Given $f \in A^{O_n}_{S^2 E}$, one has that $f =
p(a_1, \ldots, a_n)/ d^r$, where we can assume that $p(a_1,
\ldots, a_n)$ is not divisible by $d$. However, if $r > 0$ then
$p(a_1, \ldots, a_n)$ must vanish on all the diagonal matrices
with repeated eigenvalues, then $p(a_1, \ldots, a_n)$ is a
multiple of $d$, which is impossible. Therefore, $f = p(a_1,
\ldots, a_n) \in A^{Gl_n}_{{\rm End}_k (E)}$.

Let $f$ be the composite morphism ${\rm End}_k (E) \to S^2 {\rm E}
\to {\rm End}_k (E)$, $T \mapsto T \circ T^2 \circ T^t \mapsto T \circ T^2 \circ T^t \circ T_2$. The invariant functions of ${\rm End}_k (E)$  by the action by conjugation of the linear group coincide, via $f^*$, with the functions of ${\rm End}_k (E)$ that are left and right invariant by the action of the orthogonal group. The image of the morphism $k[S_m] \to {\rm End}_k (E^{\otimes m}) ={\rm End}_k (E)^{\otimes m}$, $\sigma \mapsto \tilde{\sigma}$, where $\tilde{\sigma} (e_1 \otimes \ldots \otimes e_m) = e_{\sigma(1)} \otimes \ldots \otimes e_{\sigma(m)}$, is $(({\rm End}_k (E))^{\otimes m})^{Gl_n}$, then $(S^m{\rm End}_k (E))^{Gl_n} = \langle \bar{\tilde{\sigma}} \rangle_{\sigma \in S_m}$. Now, $({\rm End}_k (E))^* = {\rm End}_k (E)$ (via the obvious isomorphism $E \otimes E^* = E^* \otimes E$), and one has in the same way that $((S^m{\rm End}_k (E))^*)^{Gl_n}$ $= \langle \bar{\tilde{\sigma}} \rangle_{\sigma \in S_m}$.

If $\alpha_1, \ldots, \alpha_n$ are the eigenvalues of $T \in {\rm
End}_k(E)$ and $\sigma = \sigma_1 \cdot \ldots \cdot \sigma_r \in
S_m$ is the decomposition of $\sigma$ as a product of disjoint
cycles (the order of $\sigma_i = m_i \geq 1$ and $m_1 + \ldots +
m_r=m$), then it is easy to check that $$ \bar{\tilde{\sigma}}(T) =
(\alpha_1^{m_1} + \ldots + \alpha_n^{m_1}) \cdot \ldots \cdot (
\alpha_1^{m_r} + \ldots + \alpha_n^{m_r} ).$$ By Newton-Girard
formulae $A_{M_n(k)}^{Gl_n} = k[a_1, \ldots, a_n] = k[S_1, \ldots,
S_n]$, where $S_i$ are the power sums symmetric polynomials (of the
eigenvalues of a matrix). Hence, a basis of $((S^m {\rm
End}_k(E))^*)^{Gl_n}$ is $\{ \bar{\tilde{\sigma}} \}_{[\sigma] \in
S'_m}$, where $S'_m = \{ [\sigma] \in S_m/\sim : \sigma $ is a product of disjoint cycles of order less than or equal to $n\}$ and $\sim$ is the conjugation relation.

Finally, $$\bar{\tilde{\sigma}} \overset{f^*}{\mapsto} \overline{T_2 \otimes \overset{m}{\ldots} \otimes T_2 \otimes T_{1 \sigma(1)}^2 \otimes \overset{m}{\ldots} \otimes T^2_{m \sigma(m)}} =: a_{\sigma} $$ where $T_{1 \sigma(1)}^2 \otimes \overset{m}{\ldots} \otimes
T^2_{m \sigma(m)} := \tau (T^2 \otimes \ldots \otimes T^2)$ and
$\tau \in S_{2m}$ is the permutation $\tau(2i-1) = 2i -1$ and
$\tau(2i) = 2 \sigma(i)$.

So we have calculated the forms $\tilde{w} \in A^*_{M_n}$ that are
left and right $O_n$-invariant. To compute the invariant integral
$w_{O_n}$ on $O_n$ it remains to impose that $\tilde{w}(I) = 0$,
where $I$ is the ideal of functions of $M_n$ vanishing on
$O_n$. Now, since $w_{O_n} \cdot \tilde{w} \cdot w_{O_n} =
\tilde{w}$, one has that $\tilde{w}(I) = \tilde{w} ( w_{O_n} \cdot
I \cdot w_{O_n})$ and $w_{O_n} \cdot I \cdot w_{O_n}$ are the
functions of $M_n$ that are left and right $O_n$-invariant and
vanish on $O_n$. These ones identify, via $f^*$, with the ideal
$I'$ of the functions of $M_n$ that are invariant by the action by
conjugation of the linear group and vanish on $Id \in M_n$.

One has $A_{M_n}^{Gl_n} = \oplus_{m \in \mathbb N} < \bar{\tilde{\sigma}} >_{[\sigma] \in S'_m} \subset A_{M_n} $ and $w_{O_n} \cdot I \cdot w_{O_n}$ identifies with $$ I' = \underset{m}{\oplus} < 1 - \frac{\bar{\tilde{\sigma}}}{\bar{\tilde{\sigma}}(Id_1)}
>_{[\sigma]\in S'_m} $$ where $Id_1 = Id \otimes
\overset{m}{\ldots} \otimes Id$. If $\sigma$ is a product of $r$
disjoint cycles (included the cycles of order $1$) it is easy to
check that $\bar{\tilde{\sigma}} (Id_1) = n^r$.

Let us denote $w_{\sigma}:= \overline{T_2 \otimes
\overset{m}{\ldots} \otimes T_2 \otimes T_{1 \sigma(1)}^2 \otimes
\overset{m}{\ldots} \otimes T^2_{m \sigma(m)}} \in S^{2m} {\rm
End}_k (E)$. In order to find $w_{O_n} = 1 + \underset{m>0}{\sum}
\sum_{[\sigma] \in S'_m} \lambda_{[\sigma]} \cdot
\dfrac{w_{\sigma}}{\bar{\tilde{\sigma}} (Id_1)}$ satisfying
$w_{O_n}(I)=0$, we have to solve, for every $m$, the system
of equations (varying $[\sigma'] \in S'_m $)
$$\left( \sum_{[\sigma] \in S'_m} \lambda_{[\sigma]} \cdot
\dfrac{w_{\sigma}}{\bar{\tilde{\sigma}} (Id_1)} \right) \left(
\frac{a_{\sigma'}} {\bar{\tilde{\sigma}}'(Id_1)} \right) = 1 .$$
If we denote $\lambda_{\sigma \sigma'} = \dfrac{w_{\sigma}
(a_{\sigma'})} { \bar{\tilde{\sigma}}(Id_1) \cdot
\bar{\tilde{\sigma}}'(Id_1)}$, then
$$(\lambda_{[\sigma]})_{[\sigma] \in S'_m} = (\lambda_{\sigma
\sigma'})^{-1}_{[\sigma],[\sigma'] \in S'_m} \cdot (1, \ldots,
1) .$$

Let us give the first three terms of $w_{O_n}$:

$$\aligned w_{O_n} & =1 + \frac{\overline{T_2 \otimes T^2}}{n}
\\ & +  \frac{(3n^2+3n+3) \cdot \overline{T_2 \otimes T_2 \otimes T^2 \otimes T^2} + (-3n-6) \cdot \overline{ T_2 \otimes T_2 \otimes T_{12}^2
\otimes T_{21}^2}} {n^4+n^3+n^2-3n} \\ & + \ldots \endaligned$$

\vskip1pt {\bf Invariant integral on $Sp_{2n}$.} \vskip3pt

Let $H_2$ be a non-singular skew-symmetric metric on a vector space
$E$ of dimension $2n$. Let $Sp_{2n}$ be the subgroup of the linear
group of the symmetries of $H_2$. In the algebraic variety
$\Lambda^2 {\rm E}^*$ of the skew-symmetric metrics, regardless of
the basis of $E$ chosen, we can define (up to a constant
multiplicative factor) the function $det$ that assigns to each
metric its determinant. So, we can consider the open set $\Lambda^2
{\rm E}^* - (det)_0$. The sequence of morphisms of varieties

$$\begin{matrix} 1 & \to & Sp_{2n} & \to & Gl(E) & \to & \Lambda^2 {\rm E}^* - (det)_0 & \to & 1 \\ & & & & & & & & \\ & & & & S & \mapsto & S^t \circ H_2 \circ S & & \end{matrix}$$ shows that $\Lambda^2 {\rm E}^*-(det)_0$ is the quotient variety of $Gl(E)$ by the symplectic subgroup $Sp_{2n}$ ($Sp_{2n}$ acting  on $Gl(E)$ on the left).

The functions of $Gl(E)$ invariants by $Sp_{2n}$ identify with the  functions of $\Lambda^2 {\rm E}^*-(det)_0$. Therefore, via the morphism of varieties ${\rm End}_k (E) \to \Lambda^2 {\rm E}^*$, $S \mapsto S^t \circ H_2 \circ S$, the functions of $\Lambda^2 {\rm E}^*$ identify with the functions of ${\rm End}_k (E)$ that are (right) invariant by $Sp_{2n}$. The morphism between the rings of functions $S^{\cdot} (\Lambda^2 E) \to S^{\cdot} ({\rm End}_k (E)^*)$ is expressed explicitly as follows $$ \begin{array}{ccl} S^m (\Lambda^2 E) & \longrightarrow & S^{2m} ({\rm End}_k (E)^*) \\ & & \\ s_1 \cdot \ldots \cdot s_m & \mapsto & \overline{ H_2 \otimes \overset{m}{\ldots} \otimes H_2 \otimes s_1 \otimes \ldots \otimes s_m} \end{array}$$ Equivalently, the left $Sp_{2n}$-invariant functions of the variety ${\rm End}_k ({\rm E})$ are the direct sum of the images of the morphisms

$$\begin{array}{ccl} S^m (\Lambda^2 E^*) & \longrightarrow &  S^{2m}
({\rm End}_k (E)^*) \\ & & \\ \omega_1 \cdot \ldots \cdot \omega_m &
\mapsto & \overline{\omega_1 \otimes \ldots \otimes \omega_m \otimes
H^2 \otimes \overset{m}{\ldots} \otimes H^2}
\end{array}$$ 
(we think of $S^{2m}({\rm End}_k (E)^*)$ as a quotient
of $E^{*2m} \otimes E^{2m}$). Therefore, the invariants of
$S^{2m}({\rm End}_k (E)^*)$ by the left and right action of $Sp_{2n}$ are

$$\aligned  & \langle \overline{\sigma( H_2 \otimes \overset{m}{\ldots}
\otimes H_2) \otimes \sigma' ( H^2 \otimes \overset{m}{\ldots}
\otimes H^2)} \rangle_{\sigma,\sigma' \in S_{2m}} \\ & = \langle
\overline{(H_2 \otimes \overset{m}{\ldots} \otimes H_2) \otimes
\sigma'(H^2 \otimes \overset{m}{\ldots} \otimes H^2)}
\rangle_{\sigma'\in S_{2m}} . \endaligned$$

Let $A_{Sp_{2n}}$ be the ring of functions of $Sp_{2n}$ and
let $w_{Sp_{2n}} \in A^*_{Sp_{2n}} \subset A^*_{M_{2n}} =
\underset{r}{\prod} S^r {\rm End}_k (E)$ be the invariant integral on
$Sp_{2n}$. The $r$-th component $[w_{Sp_{2n}}]_r$ of $w_{Sp_{2n}}$
is
$$\begin{array}{lll} \,[w_{Sp_{2n}}]_r = & \underset{{\sigma \in
S_{2m}}}{\sum} \lambda_{\sigma} \cdot \overline{(H_2 \otimes \overset m \cdots \otimes H_2) \otimes \sigma( H^2 \otimes \overset{m}{\ldots} \otimes H^2)}, & \text{ if } r=2m, \\ \,[w_{Sp_{2n}}]_r= & 0, & \text{ if } r=2m+1 . \end{array}$$

\begin{pro}{\rm (}\cite[Th. 4.3.3]{Goodman}{\rm )}
Let us consider the natural action of $Sp_{2n}$ on $E \otimes
\overset{r}{\ldots} \otimes E$, $g \cdot (e_1 \otimes \cdots \otimes
e_r) = g \cdot e_1 \otimes \ldots \otimes g \cdot e_r$. It holds
that:
\begin{enumerate}
\item $ (E \otimes \overset{2m+1}{\ldots} \otimes E)^{Sp_{2n}} = 0$.
\item $ (E \otimes \overset{2m}{\ldots} \otimes E)^{Sp_{2n}} = \langle
\sigma(H^2 \otimes \overset{m}{\ldots} \otimes H^2) \rangle_{\sigma
\in S_{2m}}$.
\end{enumerate}
\end{pro}

Let us consider the morphism $\Lambda^2 E \hookrightarrow {\rm
End}_k (E)$, ${H'}^2 \mapsto {H'}^2 \circ H_2$, that assigns to each
metric ${H'}^2$ the endomorphism associated to the pair of metrics
$H^2,{H'}^2$. Two skew-symmetric metrics are isometric (with regard
to $H_2$) if and only if their associated endomorphisms are
equivalent, and  an endomorphism $T$ (up to conjugation) is the
associated endomorphism of a skew-symmetric metric and $H_2$ if and
only if every elementary divisor of $T$ appears twice
(\cite{Ermolaev}). Let $C\subset {\rm End}_k (E)$ be the closed set of such endomorphisms. Then, $\Lambda^2E/Sp_{2n}=C/Gl(E)$. Let us write $$E=E'\oplus E',\quad H_2=\begin{pmatrix} 0 & -Id \\
Id & 0 \end{pmatrix} .$$ The diagram

$$\begin{array}{ccccccc} {\rm End}_k (E') & \hookrightarrow &
\Lambda^2E & \hookrightarrow & C & \hookrightarrow & {\rm End}_k (E) \\ T & \mapsto & \begin{pmatrix} 0 & T\\-T & 0 \end{pmatrix} & \mapsto & \begin{pmatrix} T & 0 \\0 & T \end{pmatrix} &  & \end{array}$$ shows that ${\rm End}_k (E')/Gl(E') =C/Gl(E) = \Lambda^2E/Sp_{2n}$. The ring of invariant functions of ${\rm End}_k(E')$ (by the action by conjugation of the linear group $Gl(E')$) is the image of the ring of invariant functions of ${\rm End}_k(E)$ (by the action by conjugation of the linear group $Gl(E)$), then the ring of invariant functions of $\Lambda^2E$ (by the action of $Sp_{2n}$) is the image of the ring of invariant functions of ${\rm End}_k(E)$ (by the action by conjugation of the linear group $Gl(E)$). Remember that $\{ \bar{\tilde{\sigma}} \}_{[\sigma] \in S'_m}$, where $S'_m = \{ [\sigma] \in S_m/\sim : \sigma $ is a product of disjoint cycles of order less than or equal to $n\}$, is a basis of $((S^m {\rm End}_k(E'))^*)^{Gl_n}$.

Let $f$ be the composite morphism ${\rm End}_k (E) \to \Lambda^2
{\rm E} \to {\rm End}_k (E)$, $T \mapsto T \circ H^2 \circ T^t \mapsto T \circ H^2 \circ T^t \circ H_2$. The invariant functions of ${\rm End}_k(E)$  by the action by conjugation of the linear group coincide, via $f^*$, with the functions of ${\rm End}_k (E)$ that are left and right invariant by the action of the symplectic group. Now we can calculate the invariant integral  of $Sp_{2n}$ as we have calculated the invariant integral on $O_n$.

Let us denote $w_{\sigma}:= \overline{H_2 \otimes
\overset{m}{\ldots} \otimes H_2 \otimes H_{1 \sigma(1)}^2 \otimes
\overset{m}{\ldots} \otimes H^2_{m \sigma(m)}} \in S^{2m} {\rm
End}_k (E)$ and $a_{\sigma}:= \overline{H_2 \otimes
\overset{m}{\ldots} \otimes H_2 \otimes H_{1 \sigma(1)}^2 \otimes
\overset{m}{\ldots} \otimes H^2_{m \sigma(m)}} \in S^{2m} {\rm
End}_k (E)^*$. Let $\lambda_{\sigma \sigma'} = \dfrac{w_{\sigma}
(a_{\sigma'})} { \bar{\tilde{\sigma}}(Id_1) \cdot
\bar{\tilde{\sigma}}'(Id_1)}$ and $(\lambda_{[\sigma]})_{[\sigma]
\in S_m'} = (\lambda_{\sigma \sigma'})^{-1}_{[\sigma],[\sigma']
\in S_m'} \cdot (1, \ldots, 1)$. Then, $w_{Sp_{2n}} = 1 +
\underset{m>0}{\sum} \sum_{[\sigma] \in S'_m} \lambda_{[\sigma]}
\cdot \dfrac{w_{\sigma}}{\bar{\tilde{\sigma}} (Id_1)}$.

Let us give the  first three terms of $w_{Sp_{2n}}$:
$$\aligned w_{Sp_m} & =1 + \frac{\overline{H_2 \otimes H^2}}{m}
\\ & +  \frac{(3m^2-3m+3) \cdot \overline{H_2 \otimes H_2 \otimes H^2 \otimes H^2} + (3m-6) \cdot \overline{ H_2 \otimes H_2 \otimes H_{12}^2 \otimes H_{21}^2}} {m^4-m^3+m^2+3m} \\ & + \ldots
\endaligned$$

\end{document}